\documentclass[12pt,english,reqno]{amsart}
\usepackage{palatino}
\usepackage[T1]{fontenc}
\usepackage[latin1]{inputenc}
\usepackage{pifont}
\usepackage{amssymb}
\IfFileExists{url.sty}{\usepackage{url}}
                      {\newcommand{\url}{\texttt}}

\makeatletter

 \theoremstyle{plain}
 \theoremstyle{plain}    
 \newtheorem{thm}{Theorem} 
 \theoremstyle{plain}    
 \newtheorem{lem}{Lemma} 
 \theoremstyle{remark}
 \newtheorem{rem}{Remark}
 \theoremstyle{definition}
 \newtheorem*{defn*}{Definition}
 \theoremstyle{remark}    
 \newtheorem*{claim*}{Claim}
 \newtheorem*{rem*}{Remark}

\newcommand{\PLA}{\mathrm{PLA}}
\DeclareMathOperator{\spec}{spec}

\ifx\hyperlink\undefined
\newcommand{\refp}[2]{\ref{#1}'}
\newcommand{\refpp}[1]{\ref{#1}'}
\newcommand{\hypertarget}[1]{}
\newcommand{\hrlb}{ }
\else
\newcommand{\refp}[2]{\hyperlink{#2}{\ref*{#1}'}}
\newcommand{\refpp}[1]{\ref*{#1}'}
\newcommand{\hrlb}{\\}
\fi

\AtBeginDocument{
  
}

\usepackage{babel}
\makeatother
\begin{document}

\title{Is PLA large?}

\author{Gady Kozma}

\address{GK: Institute for Advanced Study, 1 Einstein drive, Princeton NJ
08540, USA.}

\email{gady@ias.edu}
\thanks{This material is based upon work supported by the National
  Science Foundation under agreement DMS-0111298. Any opinions,
  findings and conclusions or recommendations expressed in this
  material are those of the authors and do not necessarily reflect the
  views of the National Science Foundations}

\author{Alexander Olevski\u\i}

\address{AO: Tel Aviv University, Tel Aviv 69978, Israel.}

\email{olevskii@post.tau.ac.il}

\thanks{Research supported in part by the Israel Science Foundation}

\begin{abstract}
We examine the class of functions representable by an analytic sum\begin{equation}
f(t)=\sum_{n\geq0}c(n)e^{int}\label{eq:fsumcn}\end{equation}
 converging almost everywhere. We show that it is dense but that it
is first category and has zero Wiener measure.
\end{abstract}
\maketitle

\section{Introduction}

We say that a function $f$ on the circle $\mathbb{T}$ belongs to
$\PLA$ (pointwise limits of analytic sums) if it can be decomposed
to a trigonometric series with positive frequencies (\ref{eq:fsumcn})
converging almost everywhere (a.e.)

An important observation is that the representation (\ref{eq:fsumcn})
is unique. It follows from the Abel summation and Privalov uniqueness
theorems. An analogy with the Riemannian theory suggests that the
coefficients could be recovered by Fourier formulas, provided that
$f$ is integrable (in other words, $\PLA\cap L^{1}\subset H^{1}$,
the Hardy space).

This was disproved in our note \cite{KO03}, where we constructed
a series (\ref{eq:fsumcn}) converging outside some compact $K$ of
zero measure to a bounded function $f$, but which is not $f$'s Fourier
series. Later we proved \cite{KO04,KO} that a function $f$ which
admits such a non-classic representation can be smooth, even $C^{\infty}$,
and characterized precisely the maximal possible smoothness in terms
of the rate of decrease of the Fourier coefficients.

The following density theorem is a simple consequence of these results:

\begin{thm}
\label{thm:dense}$\PLA$ is dense in the space $C(\mathbb{T})$.
Moreover, it is dense in the spaces of smooth functions $C^{k}(\mathbb{T})$
for every $k=1,2,\dotsc$ in their respective norms.
\end{thm}
The approach taken in \cite{KO03,KO04,KO} is complex-analytic, and
information is derived by examining the related function $F(z)=\sum c(n)z^{n}$
in the disk $\{|z|<1\}$. In this paper we present a purely real-analytic
construction, and use it to prove the following denseness result:
any measurable function can be carried into $\PLA$ by a uniformly
small perturbation.

\begin{thm}
\label{thm:PLA+C}Let $f\in L^{0}(\mathbb{T})$, $\epsilon>0$. Then
there is a decomposition \begin{equation}
f=g+h,\quad g\in\PLA,\,||h||_{C(\mathbb{T})}<\epsilon.\label{eq:PLA+C}\end{equation}

\end{thm}
For a stronger version, see theorem \refp{thm:PLA+C}{thm:2p} below.

This seems like a good place to compare $\PLA$ with its {}``classic''
part, the Hardy spaces. Theorems \ref{thm:dense} and \ref{thm:PLA+C}
should be contrasted against the fact that the Hardy space $H^{2}$
is closed in $L^{2}$ and has infinite co-dimension. Another interesting
fact is that $\PLA$ functions may exhibit jump discontinuities (again,
in $H^{1}$ this is impossible) --- this is a corollary from theorem
\ref{thm:PLA+C}. Finally it is worth to note that $\PLA$ contains
non-constant real functions. This comes from our approach in \cite{KO03}
where a singular inner function $I(z)$ in the unit disc was constructed
such that \[
f(t):=\frac{1}{I(e^{it})}\in\PLA.\]
Hence $f+\overline{f}=f+1/f$ gives the required example.

On the other hand we prove that $\PLA$ is rather thin in the sense
of Wiener measure and Baire category

\begin{thm}
\label{thm:category}$\PLA\cap C(\mathbb{T})$ is a set of first category
in $C(\mathbb{T})$.
\end{thm}

\begin{thm}
\label{thm:measure}The Wiener measure of $\PLA\cap C(\mathbb{T})$
is zero.
\end{thm}

\section{Small perturbations and PLA}

\subsection{Density of PLA}

First we deduce theorem \ref{thm:dense}. According to \cite{KO}
there exists a $C^{\infty}$ function $f\in\PLA$ such that\[
\widehat{f}(l)\neq0\textrm{ for some }l<0.\]
In fact the last inequality holds for infinitely many negative $l$-s.
Otherwise if $L$ is the smallest such number then \[
\sum_{n=L}^{\infty}c(n-L)e^{int}-\sum_{n=0}^{\infty}\widehat{f}(n-L)e^{int}\]
 ($c(n)$ from (\ref{eq:fsumcn})) is a non-trivial analytic sum converging
to zero a.e., contradicting Privalov's uniqueness theorem. Hence for
any $s<0$, by multiplying with an appropriate exponential we can
get an $f_{s}\in\PLA\cap C^{\infty}$ with $\widehat{f_{s}}(s)=1$.

Next, for an arbitrary $N$ consider the discrete convolution with
$e^{ist}$,\[
F_{N,s}(t):=\frac{1}{2\pi}\sum_{j=0}^{N-1}f_{s}\Big(t-2\pi\frac{j}{N}\Big)e^{is(2\pi j/N)}.\]
Since $\PLA$ is a translation invariant linear space, $F_{N,s}\in\PLA$
and\[
F_{N,s}\to f*e^{ist}=e^{ist}\quad\textrm{as }n\to\infty\]
in the $C^{k}$ norm for any $k=1,2,\dotsc$ Therefore any exponential
$e^{ist}$ admits an approximation by $\PLA$ functions, and theorem
\ref{thm:dense} follows.

\subsection{Lemmas}

We will use the technique from \cite{KO01} (where one may find additional
references and historical comments). By $L^{0}(\mathbb{T})$ we denote
the space of measurable functions $f:\mathbb{T}\to\mathbb{C}$ endowed
with the distance function\[
\rho(f,g):=\inf\{\epsilon:\mathbf{m}\{|f-g|>\epsilon\}<\epsilon\}\]
where $\mathbf{m}$ is the normalized Lebesgue measure on $\mathbb{T}$.
For a trigonometric polynomial \[
P(t)=\sum c(n)e^{int}\]
 we use the following notations

\begin{itemize}
\item $\spec P=\{ n\in\mathbb{Z}:c(n)\neq0\}$
\item $P^{*}(t)=\sup_{l<m}\left|\sum_{n=l}^{m}c(n)e^{int}\right|$.
\item $P_{[r]}(t)=P(rt)$, $r\in\mathbb{N}$.
\end{itemize}
For two trigonometric polynomials $g$ and $h$ consider the following
{}``special product'',\[
P=gh_{[r]}\textrm{ where }r>3\deg g.\]
Then the following is true (see (10) in \cite{KO01}; compare also
with \cite[sec.\ 1.1]{O1} and \cite[lemma 15]{K96}):\begin{equation}
P^{*}(t)\leq|g(t)|\cdot||h^{*}||_{\infty}+2g^{*}(t)||\widehat{h}||_{\infty},\label{eq:Pgh}\end{equation}
where $\widehat{h}$ is the Fourier transform and $||\cdot||_{p}$
is the $l^{p}$ or $L^{p}$ norm, according to context.

\begin{lem}
\label{lem:KO41}For any $\epsilon>0$ there exists a polynomial $h$
with $\spec h\subset[1,\infty[$, $\rho(h,\mathbf{1})<\epsilon$ and
$||\widehat{h}||_{\infty}<\epsilon$.
\end{lem}
The proof can be found in \cite[lemma 4.1]{KO01}.

\begin{lem}
Given a segment $I\subset\mathbb{T}$ and $\delta>0$ there is a trigonometric
polynomial $P$ such that
\begin{enumerate}
\item \label{enu:specPN}$\spec P\subset[0,\infty)$,
\item $\rho(P,\mathbf{1}_{I})<\delta$, $\mathbf{1}_{I}$ being the indicator
function, and 
\item \label{enu:P*Ieps}$P^{*}(t)<\delta$ outside of $I_{\delta}$, a
$\delta$-neighborhood of $I$.
\end{enumerate}
\end{lem}
\begin{proof}
Using lemma \ref{lem:KO41} find a polynomial $h$ such that $\rho(h,\mathbf{1})<\frac{1}{3}\delta$
and $||\widehat{h}||_{\infty}<\frac{1}{24}\delta^{2}$. Next approximate
$\mathbf{1}_{I}$ by a trigonometric polynomial $g$ such that\begin{align*}
|g(t)| & <\frac{\delta}{2||h^{*}||_{\infty}}\quad t\not\in I_{\delta} & \rho(g,\mathbf{1}_{I}) & <\frac{1}{3}\delta & ||g^{*}||_{\infty} & <\frac{6}{\delta}.\end{align*}
Finding $g$ can be done by interpolating $\mathbf{1}_{I}$ by a trapezoid
function and then taking a sufficiently large partial sum of the Fourier
expansion. Estimating $g^{*}$ can be done, for example, by noting
that a trapezoid function is a difference of two triangular function
$\textrm{T }$and for each $||T^{*}||_{\infty}\leq||\widehat{T}||_{1}=||T||_{\infty}\leq3/\delta$.
Set $P:=gh_{[r]}$. Then (\ref{eq:Pgh}) implies \ref{enu:P*Ieps}
and $\spec h\subset[1,\infty[$ implies \ref{enu:specPN}, provided
$r$ is large enough.
\end{proof}
A direct consequence:

\begin{lem}
\label{lem:Pphi}Given $\delta$ and a step function $\varphi$ which
is $0$ on a set $U$, there is a polynomial $P$ with \emph{\ref{enu:specPN}}
above such that\begin{enumerate}\setcounter{enumi}{3}\item $\rho(P,\varphi)<\delta$,\label{enu:rhoPphi}\item \label{enu:P*Uc}$\rho(P^{*},0)<\delta$
on $U$ (by which we mean $\rho(P^{*}\cdot\mathbf{1}_{U},0)<\delta$).\end{enumerate}
\end{lem}

\begin{lem}
\label{lem:PQpsi}Given $\delta$, $a$ and a step function $\psi$,
$|\psi|<a$ on $U$ there are polynomials $P$ and $Q$ such that
\emph{\ref{enu:specPN}}, \emph{\ref{enu:P*Uc}}, \begin{enumerate}\setcounter{enumi}{5}\item \label{enu:PQpsi}$\rho(P+Q,\psi)<\delta$
and \item \label{enu:Qinfa}$||Q||_{\infty}<a$\end{enumerate} are
satisfied.
\end{lem}
\begin{proof}
Fix $Q$ with \ref{enu:Qinfa} and $\rho(Q,\psi)<\frac{1}{3}\delta$
on $U$. Define a step function $\varphi$ which is $0$ on $U$ and
$\rho(\varphi,\psi-Q)<\frac{2}{3}\delta$. Now apply lemma \ref{lem:Pphi}
and get the result.
\end{proof}

\subsection{Proof of theorem \ref{thm:PLA+C}}

Let $f$ and $\epsilon>0$ be given. We may assume $\rho(f,0)<\frac{1}{4}$.
Define inductively sequences of trigonometric polynomials $\{ P_{k}\}$,
$\{ Q_{k}\}$ satisfying the conditions

\begin{enumerate}
\item \label{enu:indcta}$\rho\Big(f,\sum_{k\leq n}P_{k}+Q_{k}\Big)<4^{-n}$
\item $\spec P_{n}\subset[0,\infty[$
\item \label{enu:indctd}$||Q_{n}||_{\infty}<\epsilon2^{-n}$.
\end{enumerate}
Start with $P_{0}=Q_{0}=0$. Suppose $P_{k}$ and $Q_{k}$ are defined
with the requirements above for $k<n$. Set $f_{n}:=f-\sum_{k<n}P_{k}+Q_{k}$.
Approximate $f_{n}$ by a step-function $S_{n}$ such that \[
\rho(S_{n},f_{n})<4^{-n}.\]
We get\begin{equation}
\rho(S_{n},0)<4^{-n+2}.\label{eq:Sn0}\end{equation}
Denote $U_{n}=\{ t:|S_{n}(t)|<\epsilon2^{-n}\}$. Apply lemma \ref{lem:PQpsi}
for $\psi=S_{n}$, $U=U_{n}$, $a=\epsilon2^{-n}$ and $\delta=4^{-n}$.
We get polynomials $P_{n}$ and $Q_{n}$ such that \ref{enu:indcta}-\ref{enu:indctd}
are fulfilled for $k=n$.

Notice that if $\epsilon2^{-n}<4^{-n+1}$ then (\ref{eq:Sn0}) implies
that $\mathbf{m}U_{n}^{c}<4^{-n+2}$ and then condition \ref{enu:P*Uc}
of lemma \ref{lem:PQpsi} gives $\rho(P_{n}^{*},0)<4^{-n}$ on $U$.
These two together give \[
\rho(P_{n}^{*},0)<2^{-n}\quad\textrm{for }n\textrm{ sufficiently large.}\]
This means that the series $\sum P_{n}$ converges a.e.~and it defines
a function $g\in\PLA$. Hence denoting $h:=\sum Q_{n}$ we finish
the proof.\qed

\begin{rem}
Since theorem \ref{thm:PLA+C} is stronger than the first
part of theorem \ref{thm:dense}, one might wonder whether the second part
admits an equivalent improvement. In fact it does not, namely, there
exists a function $f\in C(\mathbb{T})$ which does not admit any decomposition
$f=g+h$ with $g\in\PLA$ and $h\in C^1$. We plan to exhibit this example in a
subsequent paper.
\end{rem}

\subsection{Representations by {}``almost analytic'' series.}

D. E. Menshov proved that any $f\in L^{0}(\mathbb{T})$ can be decomposed
to an a.e.~converging trigonometric series \begin{equation}
f(t)=\sum_{n\in\mathbb{Z}}c(n)e^{int}\label{eq:Menshov}\end{equation}
(see \cite{B64,K96,KO01}). The above technique gives the following,
{}``almost-analytic'' version:

\theoremstyle{plain}
\newtheorem*{thm2p}{Theorem \refpp{thm:PLA+C}}
\begin{thm2p}\hypertarget{thm:2p}Any $f\in L^{0}(\mathbb{T})$ can be decomposed in an almost everywhere
convergent series (\ref{eq:Menshov}) with the {}``negative'' part
$f_{-}$ converging uniformly on $\mathbb{T}$. Further, the negative
part can be taken to have arbitrarily small $U(\mathbb{T})$ norm.

\end{thm2p}We remind that the $U$-norm of a function $F$ is defined
by \[
||F||_{U(\mathbb{T})}:=\sup_{N\geq0}\Big\Vert\sum_{n=-N}^{N}\widehat{F}(n)e^{int}\Big\Vert_{\infty}.\]
We shall only sketch the proof of theorem \refp{thm:PLA+C}{thm:2p}. It requires
the technique of separating the measure error and the uniform error.
First we will need

\begin{lem}
\label{lem:KO21}For any $\gamma>0$ and $\delta>0$ there exists
a trigonometric polynomial $h$ satisfying
\begin{itemize}
\item $\widehat{h}(0)=0$, $||\widehat{h}||_{\infty}<\delta$
\item \textbf{$\mathbf{m}\{ t:|h(t)-1|>\delta\}<C\gamma$}
\item $||h^{*}||_{\infty}\leq1/\gamma$.
\end{itemize}
\end{lem}
This is lemma 2.1 from \cite{KO01}, which is the {}``non-analytic
counterpart'' of lemma \ref{lem:KO41}. Next we state a replacement
for lemma \ref{lem:PQpsi}:\theoremstyle{plain}
\newtheorem*{lem3p}{Lemma \refpp{lem:PQpsi}}
\begin{lem3p}\hypertarget{lem:4p}
Given $\gamma$, $\delta$, $a$ and a step function $\psi$, $|\psi|<a$
on $U$ there are polynomials $P$ and $Q$ such that \emph{\ref{enu:specPN}},
\emph{\ref{enu:P*Uc},}

\begin{enumerate}\item[\refpp{enu:PQpsi}]
$\mathbf{m}\{ t:|P+Q-\psi|>\delta\}<C\gamma$
and \item[\refpp{enu:Qinfa}]$||Q||_{U}<a/\gamma$\end{enumerate}
\end{lem3p}
The construction of $P$ and $Q$ is generally similar, but one has to
take $Q_{\textrm{lemma \refp{lem:PQpsi}{lem:4p}}}=
Q_{\textrm{lemma \ref{lem:PQpsi}}}h_{[r]}$
where $h$ comes from lemma \ref{lem:KO21} with the same $\gamma$
and a sufficiently small $\delta$. Notice that there is no price
to pay in lemma \ref{lem:KO21} for decreasing $\delta$.  The proof
of theorem \ref{thm:PLA+C} then applies mutatis mutandis.

\begin{rem}
Notice that if one replaces convergence almost everywhere by convergence
in measure then any function $f\in L^{0}$ admits an analytic representation
(\ref{eq:fsumcn}). This was proved in \cite{KO01}. However, here
the representation is not unique.
\end{rem}

\section{Category and measure}

Here we prove theorems \ref{thm:category} and \ref{thm:measure}.

\subsection{Relatives.}

We will use 

\begin{defn*}
(see \cite{O2}) For two functions $f$ and $g$ we say that $g$
is a relative of $f$ if there is a compact $K$ of positive measure
on the circle, and an absolutely continuous homeomorphism $h:\mathbb{T}\to\mathbb{T}$
such that\[
g(t)=f(h(t))\quad\forall t\in K.\]

\end{defn*}
In this paper the notion of relatives will be used through the following
lemma. Denote $C_{A}:=H^{\infty}\cap C$ i.e.~the space of continuous
boundary values of analytic functions on the disk.

\begin{lem}
If $f\in\PLA$ then it has a relative $g\in C_{A}$.
\end{lem}
Indeed, let $f=\sum_{n\geq0}c(n)e^{int}$. Consider the analytic extension\[
F(z)=\sum_{n\geq0}c(n)z^{n}\quad z\in\mathbb{D}.\]
According to Abel's summation theorem $F$ has non-tangential bou\-ndary
values equal to $f(t)$ at the point $z=e^{it}$ for almost every
$t$. Fix a compact $K$, $\mathbf{m}K>0$, where this limit is uniform.
Consider the so-called Privalov domain $P=P_{K}$ i.e.~the subset
of $\mathbb{D}$ created by removing, for every arc $I$ from the
complement of $K$, a disk $D_{K}$ orthogonal to $\partial\mathbb{D}$
at the end points of $I$. If $I$ is larger than a half circle, remove
$\mathbb{D}\setminus D_{I}$ instead of $D_{I}$ so in both cases
you remove the component containing $I$. Let $H$ be the Riemannian
mapping of the closed unit disc onto $P$. It is well known that $H$
belongs to $C_{A}$ and its boundary values define an absolutely continuous
homeomorphism of $\mathbb{T}$ onto $\partial P$. It is easy to see
that $g(t):=F(H(e^{it}))$ is a relative of $f$.\qed

So, theorem \ref{thm:category} will follow from

\begin{lem}
\label{lem:relatives}The set of functions with relatives in $C_{A}$
is first category in $C(\mathbb{T})$.
\end{lem}
For the proof of this lemma we need

\begin{lem}
\label{lem:epsdelM}Given numbers $\delta$, $M>0$ one can define
a positive $\epsilon(\delta,M)$ such that for any function $g$ in
$L^{2}(\mathbb{T})$, satisfying \begin{align*}
||g||_{2} & <M & |g(t)-c_{1}| & <\epsilon\textrm{ on }E_{1}\\
|c_{1}-c_{2}| & >\delta & |g(t)-c_{2}| & <\epsilon\textrm{ on }E_{2}\end{align*}
where $E_{j}$ are two disjoint sets of measure $>\delta$ and $c_{j}$
are two constants, one gets that $g$ is not in $H^{2}$.
\end{lem}
\begin{proof}
This is basically a consequence of Jensen's inequality. Assume $g$
admits an extension inside the disc as an $H^{2}$ function. Consider
the subharmonic function $G(z):=\log|g(z)-c_{1}|$. Then\begin{align*}
G(0) & \leq\int_{\partial\mathbb{D}}G=\int_{E_{1}}+\int_{\partial\mathbb{D}\setminus E_{1}}\leq\delta\log\epsilon+\int_{\partial\mathbb{D}}|g(z)-c_{1}|\leq\\
 & \leq\delta\log\epsilon+M+\frac{M}{\sqrt{\delta}}+\epsilon.\end{align*}
so taking $\epsilon$ sufficiently small we would get $|g(0)-c_{1}|<\frac{1}{2}\delta$.
Repeating the same argument for $|g(0)-c_{2}|$ leads to a contradiction.
\end{proof}

\begin{proof}
[Proof of lemma \ref{lem:relatives}]Denote (for small positive $s$
and $r$, and a large $M$) by $\mathcal{F}(s,r,M)$ the class of
all functions $f\in C(\mathbb{T})$ satisfying that there exists
\begin{itemize}
\item A $h\in\hom(\mathbb{T})$ such that for any measurable $A\subset\mathbb{T}$,
$\mathbf{m}A<s$ implies $\mathbf{m}h(A)<r/2$,
\item A compact $K\subset\mathbb{T}$, $\mathbf{m}K>r$,
\item A $g\in H^{\infty}$, $||g||_{2}<M$
\end{itemize}
such that $f\circ h^{-1}=g$ on $K$.
\begin{claim*}
$\mathcal{F}=\mathcal{F}(s,r,M)$ is nowhere dense in $C(\mathbb{T})$.
\end{claim*}
Take $f\in C(\mathbb{T})$ and a number $a>0$. Approximate $f$ uniformly
(with error $<a$) by a step-function $q$. Denote by $n$ the number
of intervals $I$ of constancy of $q$. We may assume that $1/n<s$
and that all the numbers $\{ q(I)\}$ are different. Denote\[
\delta:=\min\Big\{\frac{r}{2n},\min_{I\neq I'}|q(I)-q(I')|\Big\}\]
Find $\epsilon(\delta,M)$ from lemma \ref{lem:epsdelM}. The claim
will be proved once we show that any $u$ with $||u-q||_{2}<\epsilon$
is not in $\mathcal{F}$.

Indeed, suppose there are $h$, $K$ and $g$ as in the definition
of $\mathcal{F}$. One can see easily that for at least two intervals
$I=I_{1},I_{2}$\[
\mathbf{m}(h(I)\cap K)>\delta.\]
So, the function $g=u\circ h^{-1}$ satisfies all conditions of lemma
\ref{lem:epsdelM}, therefore it can not belong to $H^{\infty}$.
This contradiction proves the claim, and lemma \ref{lem:relatives},
and theorem \ref{thm:category} follow.
\end{proof}
\begin{rem}
The following conjecture (which would emphasize sharpness of Theorems
\ref{thm:PLA+C},\refp{thm:PLA+C}{thm:2p}) looks likely: one can construct
a function $f\in C(\mathbb{T})$ such that all its relatives $g$
satisfy the condition that $\{\widehat{g}(n)\}\not\in l^{1}(\mathbb{Z}_{-})$
or even stronger, $\not\in l^{p}(\mathbb{Z}_{-})$, $p<2$. Notice
that for $l_{p}(\mathbb{Z})$ this is true, see \cite{O2}.
\end{rem}

\subsection{Proof of theorem \ref{thm:measure}}

This follows quite easily from the Fourier representation of Brownian
motion (see e.g.~\cite[\S 16.3]{K85}). For simplicity, we will prove
it for the complex-valued Brownian bridge, i.e.~complex Brownian
motion $W$ on $[0,2\pi]$ conditioned to satisfy $W(0)=W(2\pi)$.
As is well known, $n\widehat{W}(n)$ are independent standard Gaussians
(one may take this as the definition of the Brownian bridge). We will
only use the fact that $\widehat{W}(-1)$ is independent from the
other variables. In other words we can write (as measure spaces) $C(\mathbb{T})=\Omega\times\mathbb{C}$
where $\Omega$ is the space of functions satisfying $\widehat{f}(-1)=0$
equipped with the measure of Brownian bridge conditioned to satisfy
this, and $\mathbb{C}$ equipped with the Gaussian measure. Since
$e^{-it}$ is not in $\PLA$, for any $f\in\Omega$ there can be no
more than one value $x(f)\in\mathbb{C}$ for which $f+xe^{-it}\in\PLA$.
Since the measure of a single point is always zero, we get by Fubini's
theorem that the measure of $\PLA$ is zero.

Fubini's theorem requires that $\PLA\cap C(\mathbb{T})$ be measurable.
In fact it is Borel. To show this, endow $\PLA$ with the distance
function $d(f,g)=\rho((f-g)^{*},0)$ where $f^{*}$ is defined as
for polynomials, i.e.\[
f^{*}(t)=\sup_{l<m}\left|\sum_{n=l}^{m}c(n)e^{int}\right|\quad f(t)=\sum_{n=0}^{\infty}c(n)e^{int}.\]
It is easy to see that $d$ makes $\PLA$ into a separable complete
metric space. Souslin's theorem \cite[\S 39, IV]{K66} states that
a one-to-one continuous map of a (Borel set in a) Polish space is
Borel. Using this for the identity map from $\PLA$ to $L^{0}(\mathbb{T})$
shows that $\PLA$ is a Borel subset of $L^{0}(\mathbb{T})$. The
restriction of a Borel set in $L^{0}$ to $C$ is Borel so we get
that $\PLA\cap C$ is Borel in $C$ and theorem \ref{thm:measure}
is proved.

We remark that the same argument can be used to prove theorem \ref{thm:category}.
One needs only use a version of Fubini's theorem for categories, see
\cite[theorem 15.4]{O80}.

\end{document}